\documentclass[12pt]{article}
\usepackage{amsmath,latexsym,amssymb}
\usepackage{graphicx,subfigure,amsfonts,color}
\newcommand{\torus}{{\mathbb T}}

\newtheorem{theorem}{Theorem}[section]
\newtheorem{definition}{Definition}[section]

\newcommand{\R}{\mathbb{R}}
\newcommand{\Z}{\mathbb{Z}}
\newcommand{\T}{\mathbb{T}}
\newcommand{\N}{\mathbb{N}}
\DeclareMathOperator{\sgn}{sgn}

\newcommand{\eps}{\varepsilon}

\begin{document}
\title{Weak periodic solutions and numerical case studies of the Fornberg-Whitham equation
}
\author{
G\"unther H\"ormann\thanks{Faculty of Mathematics, University of 
Vienna, A-1090 Wien, Austria} ~~\&~~
Hisashi  Okamoto\thanks{
Dept. of Math., Gakushuin University, Tokyo, 171-8588. 
Partially supported by JSPS Kakenhi 18H01137.
The present work was initiated during HO's stay in the Erwin Schr\"odinger
 Institute in Vienna, Austria.
Its support is highly appreciated. He also acknowledges the support of 
JSPS A3 foresight program: Modeling and Simulation of Hierarchical and
 Heterogeneous Flow Systems with Applications to Materials Science.  }
}

\maketitle

\begin{abstract}
Spatially periodic solutions of the Fornberg-Whitham equation are studied  to
 illustrate the mechanism of wave breaking and the formation of shocks for 
a large class  of initial data.
 We show that these solutions can be considered to be weak  solutions 
satisfying the entropy condition. By numerical experiments, we show
that the breaking waves become shock-wave type in the time evolution. 
\end{abstract}

\section{Introduction}

We denote by $\T = \R / \Z$ the one-dimensional torus group and identify 
functions on $\T$ with $1$-periodic functions on $\R$. The Fornberg-Whitham
 equation for the wave height $u \colon \T \times [0,\infty[ \to \R$ as a
 function of a spatial variable $x \in \T$ and time $t \geq 0$ reads 
\begin{equation}
u_t + uu_x + \left( 1 - \partial_x^2 \right)^{-1} u_x = 0 \qquad ( x \in \T, t > 0)
\label{eq:fw01}
\end{equation}
and is supplied with the initial condition
\begin{equation}\label{eq:fw02}
u(x,0) = u_0(x) \qquad (x \in \T).
\end{equation}

Well-posedness results, local in time, with strong solutions in Sobolev and Besov
 spaces have been obtained for both cases, periodic and non-periodic, 
in \cite{Holmes16, HolTho17}. We cannot always expect globally defined strong solutions.
This was predicted with a sketch of proof by \cite{seliger} and later proved rigorously in 
\cite{CE, Haziot17, hoermann}. These results proved the existence of wave-breaking, 
i.e., that the solution $u$ remains bounded but $u_x$ becomes singular,
 if the initial data displays sufficient asymmetry in its slope. 
However two questions seem to have remained unanswered: 
(i) What can we say about the solution  after a singularity has emerged? 
(ii) What is the nature of the singularity?  

Our first  goal is to prove existence of  globally defined 
weak solutions which extend beyond the time of wave-breaking and singularity formation. 
The second  goal of the present paper is  to show numerically that the singularity
 developing in such a solution looks very similar to a shock-wave solution of
 the inviscid Burgers equation. We emphasize that for the case of the real 
line in place of the torus, similar investigations have been carried out in \cite{FS},
 where the Fornberg-Whitham equation was formulated as and named 
Burgers-Poisson system (hence we have been unaware of that paper until 
almost completion of our current paper).

\section{Weak entropy solutions}

A well-known machinery exists for weak solution concepts for nonlinear partial 
differential equations in the form of hyperbolic conservation laws, 
but our Equation \eqref{eq:fw01} involves a non-local term and the question is
 whether or not this can be harmful to global existence. As we will prove, this is not
a big hurdle.

In the literature on nonlinear conservation laws, one of the fundamental 
references is Kru\v{z}kov's classic paper \cite{K}, where 
he considered equations of the form 
$$ 
  u_t + \left( \varphi(u) \right)_x + \psi(u) = 0 \qquad (x \in \R, t > 0).
$$
He proved, under  mild assumptions on the functions $\varphi$, $\psi$
and on the initial data $u_0$,  that a weak solution exists globally in time
and that it is unique in an appropriate class of functions. 
His setting is different from ours in two respects:
First, we have $x$ belonging to the one-dimensional torus, while Kru\v{z}kov
described the case with $x$ on the real line. Second, in Kru\v{z}kov's paper
 $\phi$ and $\psi$ are ordinary functions, while our equation involves  a 
non-local dependence on $u$ in $\psi$. The first difference causes no problem.
 The second issue indeed forces us to alter some of the technicalities along the way, 
but a close examination of the very lucid presentation of the proofs 
in \cite{K} reveals that the essence of most arguments can be applied to
 our equation with only slight modifications, which we will indicate  in the sequel. 

Note first that we may write $(1 - \partial_x^2)^{-1} u_x = K \ast u_x = K' \ast u$, 
where the convolution is in the $x$-variable only and with kernel function given 
by $K \colon \T \to \R$ is given by $K(x) = (e^x + e^{1-x})/(2 (e-1)) = \frac{\sqrt{e}}{e-1}
 \cosh(x - \frac{1}{2}) $ for 
 $0 \leq x \leq 1$ (see, e.g.,  \cite[Section 3]{hoermann}). Note that 
$K$ is continuous but is not $C^1$ on $\mathbb{T}$.
Note also that the derivative $K'$ is not continuous but is bounded. 

We will relate to Kru\v{z}kov's notation in \cite{K} as closely as possible, 
but  will switch the function arguments from $(t,x)$ to $(x,t)$ and occasionally 
write $u(t)$ to denote the function $x \mapsto u(x,t)$ with a frozen $t$. 
Compared with the main terms of the equation as labeled by Kru\v{z}kov we set 
$$
   \varphi(u(x,t)) = u(x,t)^2/2
$$ 
and 
$$
   (\psi u)(x,t) = \big((1 - \partial_x^2)^{-1} u_x(.,t)\big)(x) = (K' \ast u(.,t))(x). 
$$   
The term involving $\varphi$ conforms perfectly with the specifications from \cite{K} 
and requires no extra consideration at all. 
For the linear, but non-local, term given by the operator $\psi$ we will describe 
below suitable adaptations in the proof of uniqueness and make note of 
an alternative a priori estimate in the proof of existence. Overall in our
 situation, spatial periodicity, i.e., the compactness of $\T$, simplifies several
 estimates along the way in following the various proofs of key results in \cite{K}. 
Moreover, we have $L^\infty(\T) \hookrightarrow L^1(\T)$.

We use the convolution representation $K' \ast u$ in place of Kru\v{z}kov's 
term $\psi(u)$ in the following definition of the weak solution concept (where we also 
incorporate the initial condition into the basic inequality). 
\begin{definition}
Let $u_0 \in L^\infty(\T)$ and $T > 0$. 
  A function $u \in L^\infty(\T \times [0,T])$ is called a weak entropy 
solution of \eqref{eq:fw01}--\eqref{eq:fw02} if the following \eqref{eqn:entropy-sol} holds true:
  \begin{align}
    0 &\leq \int_{0}^{T} \int_{\T} \bigg(
|u(x,t) - \lambda| \partial_t \phi(x,t) 
+ \sgn( u(x,t)-\lambda) \frac{u^2(x,t) - \lambda^2}{2}\partial_x \phi(x,t)\nonumber\\
    &-\sgn(u(x,t)-\lambda)K'*(u(\cdot,t)-\lambda) (x) \phi(x,t) 
 \bigg) \,dx \,dt 
+ \int_{\T} |u_0(x) - \lambda|\phi(x,0)\,dx
    \label{eqn:entropy-sol}
  \end{align}  
 for any $\lambda \in \R$ and for any nonnegative $C^1$-function $\phi$ of compact support in $\T \times \R$. 
\end{definition} 

As is well-known, upon putting $\lambda = \pm \sup |u(x,t)|$ we may deduce that a weak
 entropy solution is also a weak (distributional) solution in the sense that in
 the integro-differential equation \eqref{eq:fw01} the term $u u_x$ may 
be interpreted as $\partial_x (u^2)/2$, since $u \in L^\infty$. Thus, we have 
$$
    \text{\ div}_{(x,t)} \begin{pmatrix} u \\ u^2/2 \end{pmatrix} =\partial_t u 
+ \partial_x \big(\frac{u^2}{2}\big) = - K' \ast u + u_0 \otimes \delta 
    \qquad \text{in } \mathcal{D}'(\T \times ]-T,T[),
$$
if we extend $u$ by setting it to $0$ in $\T \times ]-T,0[$. We may therefore
 call on Lemma 1.3.3 and the discussion of the weak solution concept in 
Section 4.3 in \cite{Dafermos}:  Upon possibly changing $u$ on a null-set 
 we may assume that $t \mapsto u(t)$ is continuous $[0,T] \to L^\infty(\T)$ 
with respect to the weak$^*$ topology on $L^\infty(\T)$ and we have, in particular,
$$
\lim_{t \downarrow 0}  \| u(t) - u_0 \|_{L^1} = 0,
$$
as required originally by Kru\v{z}kov in \cite[Definition 1]{K}.

\subsection{Uniqueness}
We first show that the solution of \eqref{eqn:entropy-sol} 
is unique. Existence will be proved later.

To show uniqueness of weak entropy solution we need an adaptation of 
\cite[Theorem 1]{K}---noting that for large $R$ we simply have $K = [0,T_0]
 \times \T$ and $S_\tau = \T$ in that statement---and of its proof to our
 situation with the non-local linear term $\psi u = K' \ast u$. 
This requires an appropriate replacement of the constant $\gamma$ in 
\cite[Equation (3.1)]{K} and an alternative argument in the course of the
 proof of \cite[Theorem 1]{K}, namely on the lines following \cite[Equation (3.12)]{K} concerning the term $I_4$ in Kru\v{z}kov's notation (defined there in \cite[Equation (3.4)]{K}), 
since in our case we cannot directly have a pointwise estimate calling on
 the mean value theorem for a classical differentiable function $\psi$. Instead, with two weak solutions $u$ and $v$ with initial values $u_0$ 
and $v_0$, respectively, we obviously have
$$
   \sup_{x \in \T} |K'\ast( u(\cdot, t) - v(\cdot, t))(x)| \leq \|K'\|_{L^\infty(\T)}
 \int_{\T} |u(y,t) - v(y,t)| dy
   = \|K'\|_{L^\infty(\T)} \| u(t) - v(t)\|_{L^1(\T)},
$$
which implies the following replacement of the next to last inequality 
on page 228 in \cite{K} (with mollifier $\delta_h$ and cut-off $\chi_\varepsilon$
 as chosen by Kru\v{z}kov) 
\begin{multline*}
  \int_0^{T_0} \int_{\T} \Big[ \big( \delta_h(t - \rho) - \delta_h(t - \tau) \big) 
\chi_\varepsilon(x,t) 
  |u(x,t) - v(x,t)|  
  \\ +  \|K'\|_{L^\infty(\T)} \chi_\varepsilon(x,t) \| u(t) - v(t)\|_{L^1(\T)}  \Big] dx dt \geq 0.
\end{multline*}
Therefore, we replace the inequality stretching in \cite{K} from 
the bottom of page 228 to the top of page 229 by 
\begin{multline*}
  \mu(\tau) := \int_{\T} |u(x,\tau) - v(x,\tau)| dx \leq\\  \int_{\T} |u(x,\rho) - v(x,\rho)| dx +
     \|K'\|_{L^\infty(\T)}  \int_\rho^\tau \int_\T \| u(t) - v(t)\|_{L^1(\T)} dx dt \\
     = \mu(\rho) +  \|K'\|_{L^\infty(\T)}  \int_\rho^\tau  \| u(t) - v(t)\|_{L^1(\T)} dt.
\end{multline*}
Then,  sending $\rho \to 0$, we arrive at
$$
    \| u(\tau) - v(\tau)\|_{L^1(\T)}  \leq \| u_0 - v_0\|_{L^1(\T)} +  
       \|K'\|_{L^\infty(\T)} \int_0^\tau \| u(t) - v(t)\|_{L^1(\T)}.
$$
(Here, the $L^1$-continuity of the weak solution is used.) 
Now Gronwall's inequality implies the following uniqueness result.
\begin{theorem} For any $T > 0$,  the weak solution
 to \eqref{eq:fw01}--\eqref{eq:fw02} is unique in  $\T \times [0,T]$.
 More precisely, if $u$ and $v$ are weak solutions 
with initial values $u_0$ and $v_0$, respectively, then for every $t \in [0,T]$,
$$
    \int_{\T} |u(x,t) - v(x,t)| dx \leq  e^{t \| K' \|_{\infty}} \int_{\T} |u_0(x) - v_0(x)| dx.
$$
\end{theorem}

\subsection{Existence}

We will establish the global existence here under the condition that
 $ u_0 \in L^{\infty}(\mathbb{T})$.   

The key idea in \cite{K} is to consider a parabolic 
regularization of \eqref{eq:fw01} in the form
\begin{equation}
u_t + u u_x + 
\left(  1 - \partial_x^2 \right)^{-1} u_x = 
\eps u_{xx}
\label{eq:vis01}
\end{equation}
with a small parameter $ \eps > 0$. Our aim is to show that, at least for a
 subsequence of $\varepsilon \to 0$, the corresponding solutions (all with 
the same initial value $u_0$) converge strongly in $L^1$ and are uniformly bounded in $L^{\infty}$.

For any $\eps > 0$, we will show that a strong solution $u$ 
to \eqref{eq:vis01} with initial value 
\begin{equation}\label{eq:vis01a}
  u(0) = u_0 \in L^2(\T)
\end{equation}
exists uniquely and globally in time. Due to the nonlocal term, we cannot 
directly call on the same references that Kru\v{z}kov uses in \cite{K} on
 page 231, but the desired result can be shown by a careful iteration of a 
standard contraction argument, which we describe in Appendix. The solution $u$ of 
 \eqref{eq:vis01} and \eqref{eq:vis01a}
  belongs to $ C([0,T]; L^2(\mathbb{T}))
\cap C(]0,T] ; H^1(\mathbb{T}))$ and, in fact, is smooth if $t > 0$.
Moreover, if $u_0 \in L^{\infty}(\mathbb{T}) \subset L^{2}(\mathbb{T})$, then 
$ t \mapsto \| u(t) \|_{L^\infty}$ is continuous on $[0,T]$.

Let $T > 0$ be arbitrary and consider the unique solution of \eqref{eq:vis01}
 and \eqref{eq:vis01a}.
For every $\varepsilon > 0$ we denote now the solution of
 \eqref{eq:vis01}--\eqref{eq:vis01a} by $u^{\eps}$.  
 We now have to find an alternative way to obtain
 Kru\v{z}kov's basic estimate (4.6) in \cite{K}, which he got from a 
direct application of the maximum principle. To arrive at our analogue 
of the basic estimate in \eqref{eq:maxprinc} below, we proceed as follows: 
Multiplying \eqref{eq:vis01} by $u^{\eps}$ and integrating over the spatial 
variable $x$ yields 
$$
  \int_\T u^{\eps} u_t^{\eps} dx + \int_\T (u^{\eps})^2 u_x^{\eps} dx +
 \int_\T (\psi u^{\eps}) u^{\eps} dx = \varepsilon 
\int_\T u^{\eps} u_{xx}^{\eps} dx.
$$
Noting that $u^{\eps} \mapsto \psi u^{\eps} = (1- \partial_x^2)^{-1} \partial_x 
u^{\eps} = K' \ast u^{\eps}$
 is a skew-symmetric linear operator with respect to the standard inner 
product on $L^2(\T)$, writing $u^{\eps} u_t^{\eps} = \partial_t  (u^{\eps})^2/2$
 and $(u^{\eps})^2 u_x^{\eps} = \partial_x((u^{\eps})^3)/3$, and 
integrating by parts on the right-hand side, yields (thanks to periodicity)
\begin{equation}
\frac{1}{2} \frac{d}{dt} \int_{\T} (u^{\eps})^2 dx = - \epsilon \int_{\T} (u_x^{\eps})^2 dx.
\label{eq:vis02}
\end{equation}
In particular, we have for every $t \in [0,T]$ the a priori estimate 
\begin{equation}
 \| u^{\eps}(t) \|_{L^2(\T)} \le \| u_0 \|_{L^2(\T)},
\label{eq:vis03}
\end{equation}
which is  independent of $\varepsilon$.

From \eqref{eq:vis03} and since $( 1 - \partial_x^2)^{-1}$ is a 
pseudodifferential operator of order $-2$, the $H^1$-norm of $v^\eps(t) := 
( 1 - \partial_x^2)^{-1} u^\eps_x(t)$
is bounded by $c \| u_0 \|_{L^2(\T)}$ for some constant $c$ independent 
of $\varepsilon$ and of $t$. Via the continuous imbedding $H^1(\T)
 \hookrightarrow L^\infty(\T)$ we thus have $\| v^\eps (t) \|_{L^{\infty}(\T)}
 \leq c' \| u_0 \|_{L^2(\T)}$ for some constant $c'$ independent of
 $\varepsilon$ and of $t$. We now interpret the original equation in the form 
$$ 
  \partial_t u^\eps (x,t) - \eps \partial_x^2 u^\eps(x,t)  + u^\eps(x,t)
 \partial_x u^\eps(x,t) = -  v^\eps(x,t)  
$$ 
and apply the theorem on page 230 in \cite{J} to obtain
$$
|u^\eps(x,t) | \le \| u^\eps(\cdot,s) \|_{L^{\infty}} + (T-s) \sup_{s\le \tau \le T, y \in   
   \mathbb{T}} | v^\eps(y,\tau)|
$$
for all $ 0 < s \le t \le  T$. 
(Here we regard the factor $u^\eps(x,t)$ in the term $ u^\eps(x,t) \partial_x u^\eps(x,t)$  as 
a coefficient to the first order derivative and also note that in our case
 the zero order coefficient vanishes. Although the text in \cite{J} 
 assumes more regularity of the coefficients and data, continuity is sufficient.) 
 We then let $ s \rightarrow 0$ to obtain
\begin{equation}\label{eq:maxprinc}
   \forall x \in \T, \forall t \in [0,T] : \quad 
   |u^\eps(x,t)| \leq \| u_0 \|_{L^{\infty}(\T)} + c'\, T \| u_0 \|_{L^2(\T)}.
\end{equation}
This inequality is our analogue of Kru\v{z}kov's basic estimate (4.6) in \cite{K}
 and we may now return to follow along his lines again more closely (note that 
we are closest to what he classifies as `Case A' in his paper on page 230 in next
 to the last paragraph) to establish a uniform modulus of $L^1$-continuity. 
In fact, \cite[Equation (4.7)]{K} for $w(x,t) := u^\eps(x+z,t) - u^\eps(x,t)$ holds 
in our case with $e_i =0$ and replacing the term $c w$ by $K' \ast w$, 
preserving the Lipschitz continuity properties noted on top of page 232
 in \cite{K} (in our case even globally on the compact torus). Lemma 5 
in \cite{K} is applicable in essentially the same way in establishing the key
 modulus of continuity estimate \cite[Equation (4.15)]{K} (with an appropriate
 proof variant taking into account the convolution term and showing
 (4.13) directly from (1.3) there). Thus, we have at least sketched how
 everything is in place and sufficient to apply Kru\v{z}kov's method in
 proving existence from compactness in $L^1$ (via boundedness and 
equicontinuity of the $L^1$-norm, \cite[Theorem 4.26]{Brezis2011}) and
 uniform $L^\infty$-bounds. In combination with the above uniqueness result, 
we have the following statement.
\begin{theorem}
Let $T > 0$ be arbitrary. For any $u_0 \in L^\infty(\T)$, there exists a unique 
weak entropy solution of \eqref{eq:fw01} and \eqref{eq:fw02} on the time interval $[0,T]$. 
\end{theorem}

\section{Numerical experiments by finite differences}

We now study numerical solutions and will observe that many of these exhibit 
singularities of shock-wave type. The solutions have been computed numerically under periodic boundary conditions and employing Godunov's  finite difference method.

We employ the following finite difference method with 
the uniform grid sizes  $ h = \Delta x,  \tau = \Delta t$.
The nonlinear term is discretized by 
$$
u^{n+1}_k = u_k^n - \frac{\tau}{h} \big( g(u_k^n, u_{k+1}^n) 
- g(u_{k-1}^n, u_k^n)  \big).
$$
where the numerical flux is defined by 
$$
g(u_{k-1}^n, u_k^n) =   \left\{ \begin{array}{ll}
f\left( u_{k-1}^n \right)    \qquad  &  \hbox{if} \quad 
u_{k-1}^n  \ge u_k^n ~~ \hbox{and} ~~ f\left( u_{k-1}^n\right) \ge 
 f\left( u_k^n\right),   \\
f\left( u_{k}^n \right)    \qquad  &  \hbox{if} \quad 
u_{k-1}^n  \ge u_k^n ~~ \hbox{and} ~~ f\left( u_{k-1}^n\right) \le
 f\left( u_k^n\right),   \\
f\left( u_{k-1}^n \right)   \qquad  &  \hbox{if} \quad 
u_{k-1}^n  \le u_k^n ~~ \hbox{and} ~~ f'\left( u_{k-1}^n\right) \ge 0,   \\
f\left( u_{k}^n \right)    \qquad  &  \hbox{if} \quad 
u_{k-1}^n  \le u_k^n ~~ \hbox{and} ~~ f'\left( u_{k-1}^n\right) \le 0,   \\
0 \qquad & \hbox{otherwise}. 
\end{array}
\right.
$$
Here $u_k^n$ is the approximation for $ u(kh,n\tau)$, and $f(u) = u^2/2$.
This is the Godunov method as explicated in \cite{CM}.

The nonlocal term is discretized by the following central difference scheme:
We put $ v = (1 - \partial_x^2)^{-1} u_x$, whence we have 
$ v - v_{xx} = u_x$. The function $v$ is determined by solving the 
linear equation
$$
v_k - \frac{ v_{k+1} -2 v_k + v_{k-1}}{h^2} = \frac{u_{k+1} - u_{k-1}}{2h},
$$
or, 
$$
( 1 + 2h^{-2}) v_k - h^{-2} ( v_{k+1}  + v_{k-1} )= \frac{u_{k+1} - u_{k-1}}{2h}.
$$

\begin{figure}[htbp]
\begin{center} \leavevmode
\includegraphics[width=6.5cm]{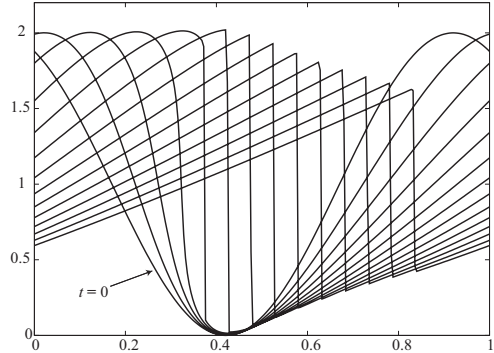}
\end{center}
\caption{ The solution from \emph{ data}1. $ 0 \le t \le  0.65$.  }
\label{zu:1}
\end{figure}

We set $ N = 1000$,  $ h = 1/N$, and  $\tau = 0.4h/q$, where $q$ is the typical 
size of the initial data. 
We first test the following two initial  data: 
\begin{align*}
& data 1  \hskip 1cm      u_0(x) = \cos (2\pi x  + 0.5) + 1,  \\
& data 2  \hskip 1cm     u_0(x) = 0.2 \cos (2\pi x) + 0.1\cos (4\pi x)   -
0.3 \sin (6\pi x)  +  0.5  
\end{align*}
The corresponding profiles of $u$ from \emph{data}1 are shown in Figure \ref{zu:1}.
As the profile moves to the right, the formation of a shock is clearly visible.
After the emergence of the shock, the jump height at the discontinuity decreases
as $t$ increases, namely, if $ \xi(t)$ denotes the position of the discontinuity at time $t$, the difference of the one-sided limits $u(\xi(t)-0) - u(\xi(t)+0)$ is a decreasing function of $t$.

\begin{figure}[htbp]
\begin{center} \leavevmode
\includegraphics[width=12.5cm]{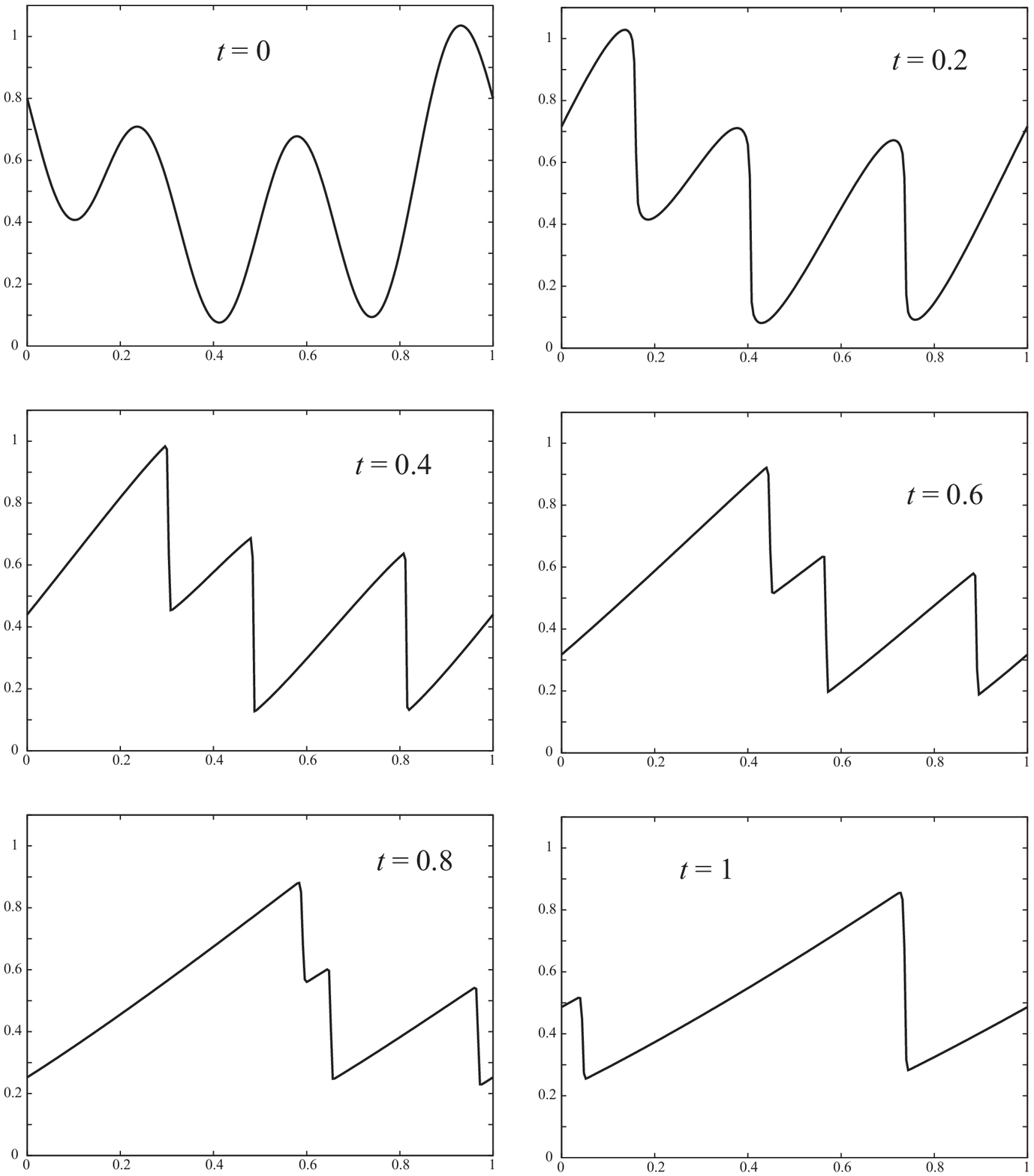}
\end{center}
\caption{ data2 }
\label{zu:2}
\end{figure}

The solution  with \emph{data}2 is shown in Figure \ref{zu:2}. In this case, even 
multiple shocks are formed and merging of some of the shocks over time can be observed.

These and other experiments suggest the development of 
wave-breaking singularities into shock discontinuities. In fact, as far as we have computed, 
only shock-type singularities in wave solutions have been found.  Moreover, the computations also suggest   that global solutions exist, if the initial data are small---the recent result in \cite[Theorem 1.5]{Itasaka} on Fornberg-Whitham-type equations with nonlinear term $\partial_x (u^p/p)$, $p \geq 5$, seems to support our observation.

\subsection{Remarks on shock conditions and asymptotic properties}

In our experiments all initial functions of the following form 
$u_0 = \sum_{1 \le k \le 3} [ a_k \cos(2k\pi x) + b_k \sin(2k\pi x) ] $
 which are periodic in $x$, but not necessarily symmetric,   produce 
a shock spontaneously, if their  $L^{\infty}$-norms are large. 
If the initial data is small, the numerical solution exists for quite a long
time. For instance, if $ u(x,0) = q\cos(2\pi x)$, a shock wave was 
observed for $ q > 0.015$. For $ q < 0.01$ the solution seems to exist
forever. For $ 0.01 < q < 0.015$, our experiments are indecisive. It
may exist forever, or it may have a shock after quite a long, 
numerically undetectable time has 
passed. Note that the guaranteed time of existence according
 to \cite[Theorem 1]{Holmes16} is inverse proportional to $\| u_0 \|_{H^2}$. 
The paper mentioned above, \cite{Itasaka} discussing a whole class 
of Fornberg-Whitham-type equations,  contains also detailed information
 on the existence time.

The propagation speed of the shock
is the same as in the case of the inviscid Burgers equation
\begin{equation}
u_t + uu_x = 0 \qquad (x \in \mathbb{T}, t > 0).
\label{eq:bur}
\end{equation}
Indeed, 
 for any $w$ in $L^2$, the function $ \left( 1 - \partial_x^2 \right)^{-1} w_x$
belongs to $H^2(\torus)$, hence is continuous, and therefore does not contribute to the Rankine-Hugoniot jump relation
$$
c = \frac{u^+ + u^-}{2}.
$$
Here, $c$ is the propagation speed of the shock and $u^+$ and $ u^-$ denote
the limit values of $u$ from the right and the left at the discontinuity,
 respectively.

Although the formation of a shock is quite similar in both equations,
the asymptotic behavior may be different. In the case of 
the Burgers equation, the jump discontinuity  vanishes in the limit and the solution 
converges to a constant function as $ t  \rightarrow \infty$ (for a proof see \cite{lax} 
or \cite[Theorems 6.4.9 and 11.12.5]{Dafermos}).
On the other hand, the Fornberg-Whitham equation possesses many traveling wave solutions 
that obviously do not decay and it is not yet conclusively shown whether
 or not for the discontinuous solutions displayed above the jump heights definitely
 decay to zero. We are not aware of a decisive result whether discontinuous
 \emph{periodic} traveling waves exist, although we discuss below a negative 
result in the case of a single shock. (Existence of non-periodic 
discontinuous wave solutions has been shown in \cite{FS,hoermann2}.)

Another feature of smooth solutions $u$ to the inviscid Burgers equation is that
 it preserves the values of maxima and minima of $ u(t,\cdot)$ as $t$ progresses. 
The Fornberg-Whitham  
equation does not keep those values constant due to the presence of the nonlocal term. 
However, our experiments suggest that, although $\max_x u(x,t)$ and $ \min_x u(x,t)$ is not 
a constant function of $t$, they are nearly constant, or, at least they seem 
to stay bounded as $t \to \infty$.

The number of peaks of any solution to the inviscid Burgers equation is 
non-increasing, while we observe that  in case of our equation they may 
increase (and decrease) as $t$ varies. 

Although several of the minor discrepancies exist as indicated, the solutions 
of \eqref{eq:fw01} have some
 similarity with solutions of the inviscid Burgers equation in
any time interval until and little after the formation of a shock discontinuity.

\subsection{Traveling waves}
In this subsection we first investigate the stability of traveling wave
 solutions and then the non-existence of such with a single shock.

Consider a continuously differentiable traveling wave $u(x,t) = U(x-ct)$ with 
speed $c > 0$. Inserting this into \eqref{eq:fw01} and integrating once we
 deduce that the profile function $U$ satisfies the equation
\begin{equation}
-cU + \frac{U^2}{2} + \left( 1 - \partial_x^2 \right)^{-1} U = 0. 
\label{eq:U}
\end{equation}
(The constant of integration may be set to zero without losing generality,
if $U$ and $c$ are suitably normalized.)

If we define $V$ by
$$
V = -cU + \frac{U^2}{2},
$$
then $2 V + c^2 = (U - c)^2 \geq 0$ and $U(x) = c \pm \sqrt{ c^2 + 2V(x)}$, 
where we choose the negative sign for the root, since only this connects
 $U(x) = 0$ with $V(x) = 0$. We obtain the differential equation 
\begin{equation}
V - V_{xx} + c - \sqrt{ c^2 + 2V} = 0,
\label{eq:V}
\end{equation}
which possesses $V \equiv 0$ as trivial solution for all $c$. 
Linearization of  \eqref{eq:V}  at $V \equiv  0$  yields
$$
V - V_{xx} - \frac{V}{c} = 0.
$$
In terms of the Fourier coefficients $(\hat{V}(m))_{m \in \Z}$ for the
 $1$-periodic function $V$, this means
$\left(1 + (2 \pi m)^2 - \frac{1}{c}\right) \cdot \hat{V}(m) = 0$ for every
 $m \in \Z$. We obtain for every $n \in \N$ a nontrivial solution to the 
linearized equation in the form $V_n(x) := \cos (2\pi n x)$, if the speed 
attains the bifurcation value
$$
   c = c_n := \frac{1}{ 1 + (2\pi n)^2},
$$
otherwise we are left with $V \equiv 0$ as the only solution.

\begin{figure}[htbp]
\begin{center} \leavevmode
\includegraphics[width=6.5cm]{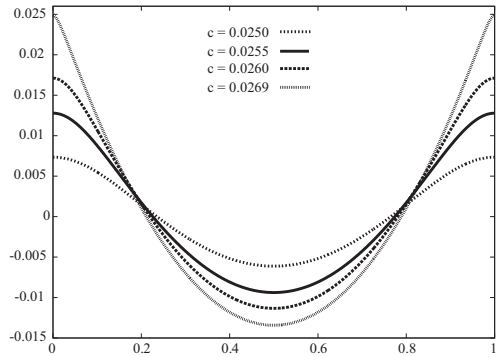}
\end{center}
\caption{ traveling wave $U$; $c = 0.025, 0.0255, 0.026, 0.0269$.  
 }
\label{zu:travel}
\end{figure}

In the case of $n=1$, nontrivial bifurcating solutions to the nonlinear
 differential equation exist in the interval $ c_1 \approx  0.0247 < c < 0.02695$. 
As $c$ tends to the upper limit, 
$\min ( c^2 + 2V)$ tends to zero, and the profile $U = c - \sqrt{c^2 + 2V}$
 tends to a function with a corner singularity. 
This Lipshitz continuous traveling wave, called peakon, has already been known 
to exist for a long time (cf.\ \cite{FW, W}) and occurs also as part of the analysis 
in \cite{ZT} . Here we have computed these traveling waves in the context above,
 i.e., as the solutions of the boundary value problem \eqref{eq:V} and illustrate 
some of them in Figure \ref{zu:travel}. 

\begin{figure}[htbp]
\begin{center} \leavevmode
\includegraphics[width=6.5cm]{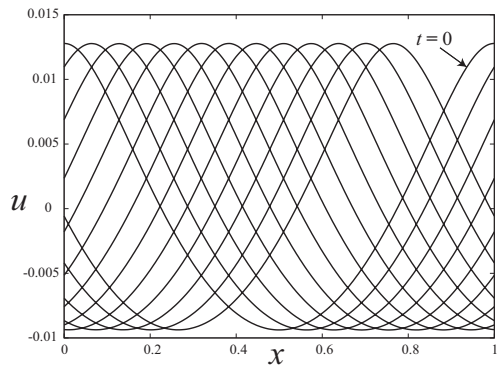}
\end{center}
\caption{ 
 The time dependent solution 
with the traveling wave as the initial data. }
\label{zu:trav}
\end{figure}

In order to investigate stability, we picked the solution corresponding to
 $ c = 0.0255$,  input it as initial data, and computed the time dependent
 solution shown in Figure \ref{zu:trav}, which illustrates that the fixed 
wave profile travels at constant speed. In the time interval used by us, 
$ 0 < t < 30$, the effect of numerical viscosity is invisible, but for long 
time intervals, it is expected to influence the numerical solution. 

\begin{figure}[htbp]
\begin{center} \leavevmode
\includegraphics[width=6.5cm]{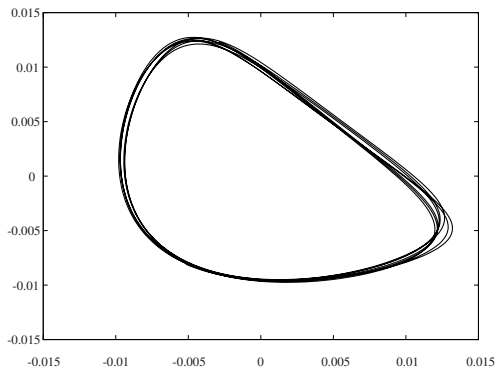}
\end{center}
\caption{  
 The time dependent solution with  $d$ as in the second case of \eqref{eq:dis2}. 
The points $( u_{300}^n,u_{600}^n)$ with $n$ corresponding to
$ 0 \le  t \le 300$ are plotted. }
\label{zu:trav4}
\end{figure}

We consider a perturbation of the initial data in the form $ U(x) + \delta d(x)$, where 
$\delta$ is $5\%$ of the amplitude of $U$ and $d(x)$ is given as follows:
\begin{equation} \label{eq:dis2}
d(x) = \cos(2k\pi x) \quad ( k = 2,3,4) \quad
\text{ or } \quad
d(x) = \left\{  \begin{array}{ll} 1 - \cos(4\pi x) \qquad  & ( 0 \le  x < 1/2), \\
0 & (\hbox{otherwise}). \end{array}  \right.
\end{equation}
Note that the latter disturbance is asymmetric. It turns out that in the time interval
 $ 0 \le  t \le 300$ the solution stays in a small neighborhood of the original
 solution $U$ and only a small oscillation could be observed as is shown 
in Figure \ref{zu:trav4}, where the points  $( u_{300}^n,u_{600}^n)$ with $n$ corresponding to
$ 0 \le  t \le  300$ were plotted. If the initial perturbation is null, 
these describe a closed curve; once the initial disturbance is added, the 
corresponding curve departs from the closed orbit, but remains in a certain neighborhood. 
These and similar experiments may be interpreted as support for the conjecture 
of stability of the traveling wave with speed $ c = 0.0255$. 
More experiments carried out for the case $ c = 0.0265$ gave similar results.

\subsubsection{Nonexistence of traveling waves with a single shock}

In view of the solutions shown earlier that created shocks after wave breaking, 
with jump height decreasing as time progresses, it is natural to ask whether
 periodic traveling wave solutions with jump discontinuities exist. For the non-periodic
 case it has been shown that discontinuous traveling waves with single jumps exist 
(see \cite{FS,hoermann2}). Contrary to this, there is no periodic traveling wave 
with a single jump, as we will argue in the following.

Suppose $U$ is the profile function for a \emph{discontinuous} traveling wave 
solution that is piecewise $C^1$ on the torus $\T$, i.e., $C^1$ except for a
 single point $x_0 \in \T$ where $U$ as well as $U'$ possess one-sided limits. 
Since the Fornberg-Whitham equation is invariant under translations in 
the $x$-variable, we may restrict to the case $x_0 = 0$, thus we may think 
of the profile function as a  $C^1$-function $U \colon [0,1] \to \R$ with $U(0) \neq U(1)$.
 The Rankine-Hugoniot shock condition requires
$$
   U(0) + U(1) = 2c.
$$

We come back to Equation \eqref{eq:U} for the profile function $U$, but 
this time keep track of the   constant of integration, for later convenience in the form
$$
   -cU + \frac{U^2}{2} + \left( 1 - \partial_x^2 \right)^{-1} U = \beta - \frac{c^2}{2} + c,
$$
with $\beta \in \R$. We now put 
$Y := U - c$ and $ Z := Y Y' = ( Y^2 /2 )'$, note that $\left( 1 - \partial_x^2 \right)^{-1} 1
 = K \ast 1 = \int_{\T} K(x) \, dx = 1$, and insert into the equation for $U$ to obtain 
$$
     \frac{Y^2}{2} + \left( 1 - \partial_x^2 \right)^{-1} Y = \beta,
$$ 
which also shows that $Z = (Y^2/2)' = - K' \ast Y$ is continuous as a function 
on the torus and piecewise $C^1$, i.e.\ can be thought of as $C^1$- function $Z \colon [0,1] \to \R$ with 
\begin{equation}
   Z(0) = Z(1). \label{eq:Zcont}
\end{equation}

We may therefore apply $(1 - \partial_x^2)$ and arrive at 
$$
   \beta = \frac{Y^2}{2} - \left( \frac{Y^2}{2} \right)'' + Y = \frac{Y^2}{2} - Z' + Y.
$$
We collect the equations for $Y$ and $Z$ in the following first-order system 
\begin{align}
  Y Y' &= Z, \label{eq:Y}\\
  Z' &= \frac{Y^2}{2} + Y - \beta = \frac{1}{2} \left( (Y + 1)^2 - (2 \beta +1) \right), \label{eq:Z}
\end{align}
for the $C^1$-functions $Y$ and $Z$ on $[0,1]$. The requirements on $U$,
 including the Rankine-Hugoniot condition, now read
\begin{equation}
    Y(0) \neq Y(1) \quad\text{and}\quad Y(0) + Y(1) = 0, \label{eq:Ycond}
\end{equation}
while for $Z$ we have the periodic continuity condition \eqref{eq:Zcont}.

We split the further analysis into two cases depending on the value of $\beta$:

\emph{Case $2 \beta + 1 \leq 0$}: Equation \eqref{eq:Z} implies that $Z' \geq 0$. 
By \eqref{eq:Zcont}, this leaves only the option that $Z' = 0$, hence
 $0 \leq (Y+1)^2 = 2 \beta + 1 \leq 0$, so that $Y$ would have to be constant
 (equal to $-1$), which is a contradiction to \eqref{eq:Ycond}.

\emph{Case $2 \beta + 1 > 0$}: As a preliminary observation, we note the following: 
The conditions \eqref{eq:Ycond} imply that either $Y(0) < 0 < Y(1)$ or $Y(0) > 0 > Y(1)$, 
hence there exists $s \in \,]0,1[$ such that $Y(s) = 0$; by \eqref{eq:Y}, 
we have also $Z(s) = 0$, thus the trajectory in the $(Y,Z)$-phase diagram passes
 through the origin $(0,0)$ and, due to \eqref{eq:Ycond} and \eqref{eq:Zcont},
 has as end point $(Y(1),Z(1)$ precisely the reflection at the $Z$-axis of
 its starting point $(Y(0),Z(0))$.

Since $2 \beta + 1 > 0$, the vector field defining the right-hand sides of
 \eqref{eq:Y} and \eqref{eq:Z} has equilibirum points in $(-1 \pm \sqrt{2 \beta + 1 },0)$
 and the qualitative analysis in \cite[Section 3]{FS} may be applied to show that no
 trajectory satisfying all the above specifications can exist (the quantities $u$, $E$, $d$ 
used there correspond here to $Y$, $Z$, $\beta$, respectively).

In summary, we have shown that there is no periodic traveling wave with a single shock. 

The question remains whether there exist traveling wave solutions with two or more 
shock discontinuities (and being piecwise $C^1$). Reviewing the above line of 
arguments, the reasoning in the first case, $2 \beta + 1 \leq 0$, seems to 
essentially stay valid (with monotonicity arguments on subintervals instead), 
whereas in the second case, $2 \beta + 1 > 0$, the crucial consequence that $Y$ 
(and hence $Z$) has to vanish somewhere is lost, if $Y$ is allowed to have 
an additional jump discontinuity. We have to leave this issue open for potential
 future analysis and note that, even in that case, $Z$ still has to be a continuous 
function on all of $\T$ thanks to the relation $Z = (Y^2/2)' = - K' \ast Y$.

\section{Concluding remarks}

The Godunov method is of first order and in further numerical studies one might
 want to employ a method of higher order such as the ENO (Essentially Non-Oscillatory)
 scheme for better accuracy. Furthermore,  the computation of 
$(1-\partial_x^2)^{-2}\partial_x u$ from $u$ in our approach may contain 
a significant truncation error. Therefore, there admittedly is a lot 
of room for improvement in the numerical experiments.
Nevertheless, we do expect our numerical solutions to correctly show several
 main qualitative features. For instance, the almost generic emergence of
 shocks in the wave solution $u$ appears to be undoubted and our experiments
 strongly indicate that many of the traveling wave solutions are stable. 
On the other hand, the large time behavior of solutions in general cannot 
be assessed quantitatively by our method and we clearly lack theoretical insight. 

In summary, we are left with (at least) the following questions which we were 
unable to answer in terms of a rigorous analysis so far:
\begin{enumerate}
\item  Existence of periodic traveling waves with jump discontinuities
 (i.e, non-decreasing shocks), although the case of a single jump could be ruled out.
\item  Boundedness of the spatial minimum and maximum of $u(t)$ as a 
function of $t$ and independent of the existence time T.
\item  Wave breaking in more generic cases than those covered by the 
asymmetry condition in the wave breaking theorems.
\item  Global (in time) existence of strong solutions for (generic) 
initial data with small Sobolev norm.
\end{enumerate}

\appendix

\section{Appendix: Proof of the existence of a global solution to the regularized equation}

We point out that Fujita and Kato's theorem on the local existence of a strong solution
of the Navier-Stokes equations (\cite{FK}) can be used to prove the existence 
of a strong solution of \eqref{eq:vis01} and \eqref{eq:vis01a}. 
Their method is explained and put into a more general context in \cite{H} and \cite{pazy}.
The same techniques have also been applied in \cite{CO}, which indeed involves
 a  nonlinearity similar to that in the Fornberg-Whitham equation. Given some familiarity with
the theory of analytic semigroups, one would quickly see that the proof described
 below is only a slight variation of the classical methods. However, we think that
 there will be some benefit or at least convenience for the reader in outlining the proof here.

We first note that $ A := - \epsilon \partial_x^2$ generates an analytic semigroup 
of operators in $L^2(\mathbb{T})$ (see, e.g., \cite{EN, kato, pazy}).  
In terms of the semigroup, the Cauchy problem  \eqref{eq:vis01} and \eqref{eq:vis01a}
 is equivalent to the following fixed point problem
\begin{equation}
u(t) = e^{-tA} u_0 + \int_0^t e^{-(t-s)A} \left[ - 
u(s)  u_x(s)  - (1-\partial_x^2)^{-1} u_x(s)  \right]  ds 
=: F(u)(t)   \label{eq:fu}
\end{equation}
We will show that this equation has a unique solution for every
$ u_0 \in L^2(\mathbb{T})$. The proof is carried out by a successive 
approximation as follows:  $u^{(0)}(t) = e^{-tA} u_0$,
 $u^{(n+1)}(t ) = F( u^{(n)})(t)$ for $ n = 0,1,\ldots$ and the solution is 
then obtained by 
showing that $F$ is a contraction mapping with respect to a suitable metric. 

In the sequel we will denote by $c$ or $c'$ various constants that may depend
 on $\epsilon$ but not on $t$.  Furthermore, let $\gamma_0$ be a constant such that
$$
\left\|  \partial_x e^{-tA} v \right\| \le \gamma_0 \| v\| t^{-1/2}  
\qquad  (v \in L^2),
$$
where here and hereafter $\| ~~\|$ denotes the $L^2$-norm.

Let $R > 0$ and suppose that we are given a continuous map $t \mapsto w(t)$,
 $[0,T] \to L^2(\torus)$, which takes its values for $t > 0$ in $H^1(\mathbb{T})$ and satisfies 
\begin{equation}
\| w(t) \| \le R \qquad ( 0 \le t \le T), 
\qquad   \| w_x(t) \| \le Rt^{-1/2}  \qquad (0 < t \le T).
\label{eq:uuu}
\end{equation}
We then have $ \| w(t) \|_{L^{\infty}} \le c \| w(t) \|^{1/2} \| \partial_x w(t) \|^{1/2}
\le c R t^{-1/4}$ and it is not difficult to deduce
\begin{align*}
\| F(w)(t)\| & \le \| u_0 \| + \int_0^t \| e^{-(t-s)A} \| 
\left[  \| w(s) \|_{L^{\infty}} \| w_x(s) \| + \| (1-\partial_x^2)^{-1} w_x(s) \|    \right] ds \\
& \le \| u_0 \| + \int_0^t  \left[ cR^2s^{-3/4}  + cR \right]  ds
= \| u_0 \| + 4 cR^2 t^{1/4} + cR t.
\end{align*}
Similarly, 
\begin{align*}
\| \partial_x  F(w)(t)\| & \le \gamma_0\| u_0 \| t^{-1/2} + \int_0^t \| \partial_x e^{-(t-s)A} \| 
\left[  \| w(s) \|_{L^{\infty}} \| w_x(s) \| + \| (1-\partial_x^2)^{-1} w_x(s) \|    \right] ds \\
& \le \gamma_0  \| u_0 \|t^{-1/2} + \int_0^t \gamma_0  
(t-s)^{-1/2} \left[ c R^2 s^{-3/4}  + cR \right]  ds \\
&
= \gamma_0  \| u_0 \| t^{-1/2} + \gamma_0 c R^2 t^{-1/4}B(1/2,1/4) +
 2 \gamma_0 c R t^{1/2},
\end{align*}
where $\gamma_0  $ is as above  and $B(\cdot,\cdot)$ denotes Euler's beta function.

Now we may take any $R > \max\{  \| u_0 \|, \gamma_0 \| u_0 \| \}$. 
Then we may choose $T$ small enough such that
 any $w$ satisfying \eqref{eq:uuu} implies that \eqref{eq:uuu} holds also with $F(w)$  in place of $w$. 

If we equip $Y_T := C([0,T] ; L^2) \cap C(]0,T] ; H^1)$ with the (complete) norm 
$$
\| u \|_* := \max \left\{ 
\max_{0 \le t \le T} \| u(t) \|, \sup_{0 < t \le T} t^{1/2} \| u_x(t) \| \right\},
$$
then the corresponding closed ball $B_R$ of radius $R$ around $0$ in $Y_T$ is 
mapped into itself (upon noting that $t \mapsto \partial_x F(u)(t)$ is 
continuous $]0,T] \to L^2$ for every $u \in Y_T$).

Our next task is to show that the map $F$ is a contraction with respect to 
the norm $\| ~~ \|_*$ (for some $T > 0$ and on bounded subsets). 
The proof is similar to the above estimates. Indeed, 
\begin{align*}
 \| F(w)(t) - F(z)(t) \| \le & \int_0^t 
\| e^{-(t-s)A} \| 
\bigg[  \| w(s) \|_{L^{\infty}} \| w_x(s)  - z_x(s) \| \\
& + 
\| w(s) - z(s) \|_{L^{\infty}} \| z_x(t) \| +   
\left\| (1-\partial_x^2)^{-1} \big( w_x(s) - z_x(s) \big)  \right\|    \bigg] ds \\
& \le \int_0^t \Big(  c\| w\|_* \| w-z \|_*  s^{-3/4} + 
  c\| z\|_* \| w-z \|_*  s^{-3/4} +  c \| w-z \|_* 
\Big) ds \\
& \le c' t^{1/4} \big(  \| w \|_* + \| z\|_* \big) \| w-z \|_* + ct \| w- z\|_*.
\end{align*}
and similarly for $ \| \partial_x F(u) - \partial_x F(z)\| $. Thus, for  
$T$ sufficiently small,  $F$ becomes a contractive mapping on $B_R$ and we
 obtain a unique solution $u \in Y_T$ of \eqref{eq:fu}.

\begin{theorem}
Let $R>0$ and $\eps > 0$, then there exists  $T >0$ such that 
 \eqref{eq:vis01} and \eqref{eq:vis01a} possesses a unique 
strong  solution  in $[0,T]$ for every $u_0 \in L^2(\torus)$ such 
that $\max\{  \| u_0 \|, \gamma_0 \| u_0 \| \} < R$.  
\end{theorem}

We emphasize that here $T$ depends on $R$ and $\epsilon$, but not on any
individual $u_0$ as far as $u_0$ satisfies $\max\{  \| u_0 \|, \gamma_0 \| u_0 \| \} < R$.
Since the spatial $L^2$-norm in the solution is conserved over time,
 the solution may be  continued any number of times, thus we obtain a solution globally in time. 
Therefore, we conclude the global unique existence of a solution to the parabolic equation 
\eqref{eq:vis01}. The following additional features of this solution are used
 in the proof of the weak solution: Due to the properties of the heat kernel,
$$
u \in C^\infty(\mathbb{T} \times\, ]0,T])
$$
and moreover, if $u_0 \in L^\infty$, although $t \mapsto u(t)$ may be 
discontinuous at $ t = 0$ as a map into $L^{\infty}$, the real-valued map 
$t \mapsto \| u(t) \|_{L^{\infty}}$ is continuous on the closed interval $[0,T]$.

\end{document}